\documentclass{my_m2an} 

\usepackage[utf8]{inputenc}
\usepackage[T1]{fontenc}
\usepackage{fullpage}
\usepackage{amsmath,amssymb}
\usepackage[mathcal]{euscript} 
\usepackage{stmaryrd}
\usepackage{tikz}
\usepackage{url}

\usepackage{color}
\definecolor{bfonce}{rgb}{0.,0.,0.8}	
\definecolor{bclair}{rgb}{0.87,0.92,1.}
\definecolor{vertf}{rgb}{0,0.55,0.1}
\definecolor{paille}{rgb}{1.,1.,.4}
\definecolor{rougec}{rgb}{1,0.4,0.}
\definecolor{or}{rgb}{0.98,0.6,.1}
\definecolor{vclair}{rgb}{0.6,0.95,0.75}
\definecolor{orangec}{rgb}{1.,0.6,0.}
\definecolor{rougef}{rgb}{0.8,0.,0.}

\usepackage[normalem]{ulem}
\newcounter{corr}
\definecolor{violet}{rgb}{0.580,0.,0.827}
\newcommand{\corr}[3]{\typeout{Warning : a correction remains in page
\thepage}
				\stepcounter{corr}        
				{\color{blue}\ifmmode\text{\,\sout{\ensuremath{#1}}\,}\else\sout{#1}\fi}
       {\color{red}#2}
       {\color{violet} #3}}

\newcommand{\bfe}{{\boldsymbol e}}

\newcommand{\bfF}{{\boldsymbol F}}
\newcommand{\bfn}{\boldsymbol n}
\newcommand{\bfu}{{\boldsymbol u}}

\newcommand{\bfx}{\boldsymbol x}

\newcommand{\edgespart}{{\mathfrak F}}
\newcommand{\dx}{\ \mathrm{d} x}
\newcommand{\dt}{\ \mathrm{d} t}

\newcommand{\edge}{\sigma}
\newcommand{\edgeb}{\tau}

\newcommand{\mesh}{{\mathcal M}}

\newcommand{\dive}{{\mathrm{div}}}

\newcommand{\gradi}{\boldsymbol \nabla}

\newcommand{\mnn}{{m\in\xN}}
\newcommand{\nnn}{{n\in\xN}}

\newcommand{\exm}{^{(m)}}

\newtheorem{theorem}{Theorem}[section]
\newtheorem{lemma}[theorem]{Lemma}

\theoremstyle{definition}

\theoremstyle{remark}

\begin{document}
\title[Analysis of finite volume schemes]{Finite volume schemes and Lax--Wendroff consistency}
 \author{R. Eymard, T. Gallou\"et, R. Herbin and J.-C. Latché} 

%
 
%
 \begin{abstract}
 
 \medskip
 
We present a  (partial) historical summary of the mathematical analysis of finite difference and finite volume methods, paying  special attention to the Lax--Richtmyer and Lax--Wendroff theorems.
We then state a Lax--Wendroff consistency result for convection operators on staggered grids (often used in fluid flow simulations), which illustrates a recent generalization of the flux consistency notion designed to cope with   general discrete functions.
\end{abstract}
\maketitle
\begin{center}
\emph{This work is dedicated to the memory of our dear colleague and friend Anton\'{\i}n Novotn\'y, \\ who suddenly passed away on June 3, 2021.}
\end{center}
\medskip
%

\medskip
 
%

%
The finite volume method (FVM) has been used for over 60 years in computational fluid mechanics, and more than 30 years in solid mechanics. 
However, while the finite element method, also introduced in mechanics, had already been the object of several mathematical works at the turn of this millenium, the American Mathematical Society (AMS)  2000 Classification only mentions the term ``finite volume''  in sections 74S10 (Mechanics of deformable solids) and 76M12 (Fluid mechanics); it was not until 2010 that this term appeared in section 65 `` Numerical analysis '' (65M08 , 65N08).

In some respects, the FVM is close to the finite difference method (FDM): for instance it has often been called ``conservative finite difference method'' by the hyperbolic numerical community, and ``finite difference method'' in the oil industry.
It might be for this reason that the first attempts to show the convergence of the FVM were to try and copy the FDM framework. 
We show in the first section below that the famous Lax-Richtmyer theorem developed for linear FD schemes fails to give an adequate answer, even for a linear finite volume (FV) scheme.
We then turn to the Lax--Wendroff theorem that gives two fundamental tools for the analysis of FV schemes, and recapitulate the main steps of the proof of convergence of the schemes that were derived in the 90's. 
In Section 4, we extend the Lax--Wendroff theorem to the case of staggered meshes and show how it can be used for the celebrated MAC grid.

\section{The finite difference method and the  Lax-Richtmyer theorem} \label{sec:laxr}

In classical  numerical analysis textbooks, we are taught that, in order to show the convergence of a finite difference (FD) scheme, one should show its stability and its consistency.
The founding result in this regard is the Lax--Richtmyer theorem \cite {lax-ritchmyer} due to P.D. Lax, who exposed it in a seminar at NYU in 1954.
The so-called Lax equivalence theorem can be summarized as follows (see for example \cite[Theorem 1.5.1]{strickwerda}):
\begin{theorem}[Lax--Richtmyer]\label{theo:lr}
 Consider a linear partial differential equation for which the initial value problem is well posed, and a finite difference scheme consistent for its approximation; then this scheme is convergent if and only if it is stable. 
\end{theorem}
 In some articles and textbooks (see for instance \cite[p. 142]{leveque}), maybe because of the name ``equivalence theorem'', and because uniform meshes are considered, the Lax--Richtmyer theorem is replaced by the equivalence
\begin{equation}\label{betise}
 \mbox{ ``consistency + stability} \Longleftrightarrow \mbox{convergence''}
\end{equation}
 However, Theorem \ref{theo:lr} does not state the equivalence \eqref{betise}, which, in fact does not hold in the general case of a PDE which is discretized with a non constant space step, even in the linear case.
As an example,  let us consider the approximation on $\xR\times]0,T[$, where $T>0$ is the final time of the study, of the linear transport equation
\begin{align}
&\partial_t u (x,t) +  \partial_x (a u) (x,t)  = 0,  \,  x \in \xR, 
\, t \in ]0,T[ \label{equ},  \\ 
&u(x,0)  = u_{\rm ini}(x), \,  x \in \xR,
\label{ci} 
\end{align}
for some given $a>0$, and initial data $u_{\rm ini} \in C^\infty_c(\xR,\xR)$, with $\partial_t$ (resp. $\partial_x$) the time (resp.  space) partial derivative. 
Let $\bar u(x,t) = u_{\rm ini}(x-at)$  be the exact (unique) solution of this problem.
The FDM applied to \eqref{ci} is classically defined by choosing a strictly increasing sequence $(x_i)_{i\in\xZ}$ of real numbers, such that $h := \max_{i\in\xZ}(x_{i+1}-x_i) < \infty$ and $\underline{h} =\min_{i\in\xZ}(x_{i+1}-x_i)>0$, and a time step $\delta t  =T/N$, for $N\in\xN$ with $N>1$ (here a constant time step is considered for simplicity). 
The initial data is discretized by defining the following quantities that depend on $h$ implicitly:
\begin{equation}\label{eq:inistandard}
 u^0_i= u_{\rm ini}(x_i)\textrm{ for any }  i \in \xZ.
\end{equation}
Since some upwinding is indeed necessary for stability purposes, the approximation of $\partial_x u$ at the point $x_i$ is upwinded, so that the FD scheme reads

\begin{equation}\label{eq:fdmstandard}
\frac  {u_i^{n+1}-u_i^n} {\delta t}+ a\frac {u_i^n- u_{i-1}^n} {h_{i-1/2}}=0 \textrm{ for any }  i \in \xZ \textrm{ and } n \in \llbracket 0, N-1 \rrbracket,
\end{equation}
with $h_{i-1/2}=x_i - x_{i-1}$.
In this context, the definition of the terms \emph{consistency}, \emph{stability} and \emph{convergence} is the following:
 \begin{itemize}
    \item \emph{Consistency}: it requires two conditions:
    \begin{itemize}
     \item the consistency of the discretization  of the initial condition, that is 
     \begin{equation}\label{eq:consistinifdm}
      \lim_{h\to 0} \max_{i \in \xZ} |u^0_i-u_{\rm ini}(x_i)|=0.
     \end{equation}
     \item the consistency of the discretization of the PDE \eqref{equ}; setting $t_n = n\delta t$ for $n \in \xN$, it reads:
    \begin{equation}\label{eq:consistfdm}
     \lim_{\substack{h\to 0 \\ \delta t\to 0}}\max_{i \in \xZ, n \in \llbracket 0, N-1 \rrbracket} \big |\frac  {u(x_i,t_{n+1})-  u(x_i,t_{n})} {\delta t }+ a\frac {u(x_i,t_{n})- u(x_{i-1}^n,t_{n})} {h_{i-1/2}}\big| = 0.
    \end{equation}
    \end{itemize}
    \item \emph{Stability}:  There exists $C$ depending only on $u_{\rm ini}$ (and thus not on $h$ nor on $\delta t$) such that
    \begin{equation}\label{eq:stabfdm}
     \max_{i \in \xZ, n \in \llbracket 0, N-1 \rrbracket} |u^n_i | \le C.
    \end{equation}
    \item \emph{Convergence}: $\max_{i \in \xZ, n \in \llbracket 0, N-1 \rrbracket} |u^n_i -u(x_i,t_{n})| \to 0$ as $h \to 0$.
\end{itemize}
Clearly, condition \eqref{eq:inistandard} ensures the consistency condition \eqref{eq:consistinifdm} on the initial condition. 
Moreover, the consistency condition \eqref{eq:consistfdm} on the PDE can be obtained  by Taylor expansions. 
Finally, the stability condition \eqref{eq:stabfdm} can be shown under the CFL condition $a \delta t\le \underline{h}$. 
The  conditions for the Lax-Richtmyer to hold are therefore satisfied, and the convergence of the scheme \eqref{eq:inistandard}-\eqref{eq:fdmstandard} is thus proven. 

Consider now a variant of the scheme \eqref{eq:fdmstandard} obtained by keeping \eqref{eq:inistandard}, but replacing $h_{i-1/2}$ by $h_i = \dfrac {x_{i+1} - x_{i-1}} 2$ in \eqref{eq:fdmstandard}:
\begin{equation}\label{eq:vfstandard}
\frac  {u_i^{n+1}-u_i^n} {\delta t}+ a\frac {u_i^n- u_{i-1}^n} {h_{i}}=0 \textrm{ for any }  i \in \xZ \textrm{ and } n \in \llbracket 0, N-1 \rrbracket,
\end{equation}

In the specific case where $x_{2k+1}-x_{2k} = h/2$ and $x_{2k+2}-x_{2k+1} = h$ for all $k\in\xZ$, we get that $h_{i} = \frac {3} {4} h$ for all $i\in\xZ$ and the consistency property no longer holds. 
Therefore, if the equivalence \eqref{betise} were true, the scheme  \eqref{eq:inistandard}-\eqref{eq:vfstandard} would not be convergent in the above sense.
However, let us show that this scheme is in fact convergent in the same sense. 
Let us write the finite difference scheme obtained at the points  $\widetilde x_i = \frac {x_{i} + x_{i+1}} 2$ instead of the points $x_i$, defining the values $\widetilde u_i^n$:

\begin{equation}\label{eq:tildeinistandard}
\widetilde  u^0_i= u_{\rm ini}(\widetilde x_i)\textrm{ for any }  i \in \xZ,
\end{equation}
and
\begin{equation}\label{eq:tildefdmstandard}
\frac  {\widetilde u_i^{n+1}-\widetilde u_i^n} {\delta t}+ a\frac {\widetilde u_i^n- \widetilde u_{i-1}^n} {h_{i}}=0 \textrm{ for any }  i \in \xZ \textrm{ and } n \in \llbracket 0, N-1 \rrbracket.
\end{equation}
Then the scheme \eqref{eq:tildeinistandard}-\eqref{eq:tildefdmstandard} again satisfies  the consistency condition \eqref{eq:consistinifdm},  and satisfies  the consistency condition \eqref{eq:consistfdm} since $h_i =\widetilde  x_i -\widetilde  x_{i-1}$.
The stability condition \eqref{eq:stabfdm} can again be shown under the CFL condition $a \delta t\le \underline{h} \le \frac {3} {4} h$. 
Hence the Lax-Richtmyer theorem  yields the convergence of the values $\widetilde u_i^n$ to $u(\widetilde  x_i,t_{n})$.
Using $|\widetilde u_i^0 - u_i^0|\le  h \max|u_{\rm ini}'|$, the maximum principle applied to the difference $\widetilde u_i^n - u_i^{n}$, solution to the equation obtained by subtracting \eqref{eq:vfstandard} to \eqref{eq:tildefdmstandard}, shows that 
\[
 \max_{i \in \xZ, n \in \llbracket 0, N-1 \rrbracket} \big |\widetilde u_i^n - u_i^{n}\big|\le h \max|u_{\rm ini}'|,
\]
which implies that
\[
 \max_{i \in \xZ, n \in \llbracket 0, N-1 \rrbracket} |u^n_i - u(x_i,t_n)|\le  2h \max|u_{\rm ini}'| + \max_{i \in \xZ, n \in \llbracket 0, N-1 \rrbracket} |\widetilde u^n_i -  u(\widetilde x_i,t_n)|,
\]
leading to the convergence of the scheme in the same sense as above.
This example shows that on the one hand the direction $\Leftarrow$ of \eqref{betise} is not true, and that on the other hand the Lax-Richtmyer theorem cannot be applied directly to obtain the convergence of such a scheme, since the scheme is not consistent in the above defined sense.  
In fact, the scheme \eqref{eq:vfstandard} is the 1D upwind FV scheme with the control volumes $]\frac {x_{i-1} + x_{i}} 2,\frac {x_{i} + x_{i+1}} 2[$, as shown in the next section.
In the scheme \eqref{eq:fdmstandard}, the partial derivative $\partial_x u$ is upwinded. 
In the scheme \eqref{eq:vfstandard}, the unknown $u$ itself is upwinded.

Note that the equivalence \eqref{betise} also does not hold in the case of an elliptic operator, see the example of a non uniform 1D mesh in \cite{matapli} or \cite[section 5.2]{book}.
For the analysis of such schemes, other notions must be introduced, and it seems that Peter D. Lax first identified them, in collaboration with B. Wendroff.
 
\section{The finite volume method and the Lax--Wendroff theorem} \label{sec:loi+vf}

In a very famous article of 1960, P.D. Lax and B. Wendroff \cite{lax-60-sys} consider discretization schemes for nonlinear hyperbolic systems of conservation laws and show that if a conservative scheme with consistent fluxes, in the sense that they define and which is stated below (see Section \ref{sec:flux-cons}), converges a.e. and boundedly towards a limit, then this limit is necessarily a weak solution of the system.
We call this property \emph{Lax--Wendroff (LW) consistency}.
The notions of \emph{flux consistency} and \emph{flux conservativity} highlighted in  \cite{lax-60-sys}  are truly fundamental for the convergence analysis of the FVM for the hyperbolic equations considered in \cite{lax-60-sys}, as are their extensions to elliptic and parabolic conservation equations.
In order to explain these terms, we consider, again in the 1D case, a general differential form of a conservation law, written on the whole space  $\xR$ and on a time interval $]0,T[$ where $0 < T < +\infty$ is the final time:
 \begin{equation}
 \partial_t u(x,t) +  \partial_x \mathbb F (x,t) = 0,
\label{conserv.loc} 
\end{equation}
stating the conservation of the quantity $u$ at each point $x\in\xR $ and each time $t\in ]0,T[$,  with $\mathbb F$ a vector function depending only on  $x$ and $t$ through the unknown $u$. 
 In addition to this equation, an initial condition is assumed to be given on $u$.
Simple examples of such a conservation law include 
\begin{itemize}
 \item the transport equation $ \mathbb  F (x,t) = a u(x,t)$ (linear hyperbolic equation, introduced in the previous section),
 \item the heat equation $ \mathbb F (x,t) = -\partial_x u(x,t)$ (linear parabolic equation),
\item Burgers' equation $ \mathbb F (x,t) =   u^2(x,t)$ (nonlinear hyperbolic equation),
\item the porous media equation $ \mathbb F (x,t) = - \partial_x u^p (x,t)$, $p>1$  (nonlinear parabolic equation).
\end{itemize}
The FVM consists in approximating the integral form of the conservation law, that is to say the balance on the time-space cuboid $] x, x + \delta x [\times] t, t + \delta t [$ (for given $\delta x>0$ and $\delta t>0$), rather than the PDE itself (this corresponds to the way such an equation is derived from physical conservation principles).
Note that with a non zero source term on the right hand-side of \eqref{conserv.loc}, the equation is then usually renamed ``balance law'', and this modification poses no problem for a FV discretization.
The integral form relative to the differential form \eqref{conserv.loc} of a conservation law reads
 \begin{equation}\label{bilan}
 \int_x^{x+\delta x} (u(x,t+\delta  t) - u(x,t)) \dx + \int_{t}^{t+\delta t} \Bigl(\mathbb  F (x+\delta x,t) - \mathbb  F (x,t)\Bigr) \dt  = 0.  
\end{equation}
 Let $(]x_{i-1/2},x_{i+1/2}[)_{i \in \xZ}$ be a family of intervals of $\xR$ (also called control volumes or grid cells), with $x_{i-1/2} < x_{i+1/2}$ and such that $\cup_{i \in \xZ}]x_{i-1/2},x_{i+1/2}[= \xR$ (or $\Omega\subset \xR$), and let $\delta t  =T/N$, $N>1$ (the time step could be taken non constant).
Let $x_{i}$ be some chosen point in the cell $]x_{i-1/2},x_{i+1/2}[$ (the choice of this point is constrained by the flux consistency property in the case of elliptic or parabolic problems, but not in the hyperbolic case, contrarily to the FD consistency in the sense of the previous section).
The discrete unknowns, $\{u_i^n,$ $i \in \xZ$, $n \in \{0, \ldots,N\}\}$ are expected to be approximations of $u(x_i,t_n)$ with $t_n=n\delta t$. 
The integral form \eqref{bilan} is written on each control volume $]x_{i-1/2},x_{i+1/2}[$ and time interval $]t_n,t_{n+1}[$, leading to the following FV scheme (with an explicit time scheme for the sake of simplicity): 

\begin{equation}\label{vf}
h_i \frac  {u_i^{n+1}-u_i^n} {\delta t}+ {F_{i+\frac 1 2}^n- F_{i-\frac 1 2}^n} =0 \textrm{ for any }  i \in \xZ \textrm{ and } n \in \llbracket 0, N-1 \rrbracket,
\end{equation}
where $h_i = x_{i+\frac 1 2}- x_{i+\frac 1 2}$, and $F_{i+\frac 1 2}^n$ is the numerical flux, which will be expressed in terms of the discrete unknowns to yield a numerical approximation of $\mathbb F (x_{i+1/2}, t_n)$.
 Note that  $F_{i+\frac 1 2}^n$ is the numerical flux outgoing $]x_{i-1/2},x_{i+1/2}[$ to the right and its opposite is  the numerical flux outgoing $]x_{i+1/2},x_{i+3/2}[$ to the left: this copies the situation of the exact flux $\mathbb F (x_{i+1/2}, t_n)$.
This is the well known ``\emph{local conservativity}'' or ``\emph{flux conservativity}'' property, which is important in physical applications, but also fundamental in the  mathematical analysis of the FVM. 
Indeed, it is thanks to this property that we may hope to prove some convergence properties of the method, both for elliptic or parabolic type and for hyperbolic equations despite the loss of the consistency in the FD sense, as presented in the first section of this paper.

\subsection{Flux conservativity}
Writing the FV scheme in one dimension naturally ensures the numerical flux conservativity, since only one flux is defined: $F_ {i + \frac 1 2} $ at the interface $x_{i + \frac 1 2} $.
In a multi-dimensional framework, ($d = 2 $ or 3) the PDE \eqref {conserv.loc} is written as $\partial_t u + \dive \mathbb F = 0 $, where $\mathbb F$ is a vector function of $x$ and $t$ and $\dive$ the space divergence operator.
The scheme \eqref{vf} is now written for a control volume $K$:
\begin{equation}\label{vf-md}
|K| \frac  {u_K^{n+1}-u_K^n} {\delta t} + \sum_{\sigma \subset \partial K} |\sigma|F_{K,\sigma}^n =0,   \;   n \in \llbracket 0, N-1 \rrbracket,
\end{equation}
where $|K|$ (resp. $|\sigma|$) is the volume or surface of $K$ (resp. the surface or length of $|\sigma|$) and  $|\sigma| F_{K,\edge}^n$ is the numerical outgoing flux from $K$ through the face $\edge$; it is an approximation of the ougoing normal flux $\int_{\edge} \mathbb  F (x,t_n) \cdot \bfn_{K,\edge}$  (where $\bfn_{K,\edge}$ is the normal vector to $\edge$ outward $K$), which is expressed in terms of the discrete unknowns $(u_M^n)_{M \in \mesh}$ :
\[
 F_{K,\edge}^n = \mathcal F_{K,\edge}^n ((u_M^n)_{M \in \mesh}).
\]
In 2 or 3D, the numerical flux is defined on either side of the  interface $\edge$.
Suppose that the interface $\edge$ separates the control volumes $K$ and $L$, which we write as $\edge=K|L$, then the flux conservativity reads
  \begin{equation}\label{conserv}
 F_{K,\sigma}^n =  - F_{L,\sigma}^n.
\end{equation}

\subsection{Flux consistency} \label{sec:flux-cons}
Let us turn back to the 1D case, for ease of notations.
The numerical flux is said to be consistent if for a sequence of time and space discretizations, indexed by $m$ and such that  $h\exm \to 0$ and $\delta t\exm \to 0$ as $m \to +\infty$, one has
\[
\max_{i \in \xZ, n \in \llbracket 0, N-1 \rrbracket}  |\mathbb F(x_{i+\frac 1 2},t_n)  - \widetilde F_{i+\frac 1 2}^n| \to 0  \mbox{ as } m \to + \infty,
\]
where $\widetilde F_{i+\frac 1 2}^n$ is the quantity obtained from $F_{i+\frac 1 2}^n$ when replacing the discrete unknowns by the values of a regular function $u$:
\[
 \widetilde F_{i+\frac 1 2}^n = \mathcal F_{i+\frac 1 2}^n\Bigl((u(x_i,t_n))_{i\in \xZ}\Bigr).
\]
In the context of nonlinear hyperbolic equations,  $\mathbb{F} (x_{i+\frac 1 2},t_n) = f(u (x_{i+\frac 1 2},t_n))$ with $f \in C(\xR,\xR)$, this definition of consistency is equivalent to the usual notion of consistency introduced by Lax (written here for a two-point scheme):
\[
\mathcal F_{i+\frac 1 2}^n (u,u) = f(u) 
\]
if one assumes $\mathcal F_{i+\frac 1 2}^n$ to be Lipschitz continuous, or at least,  ``lip-diag'' see \cite[Remark 5.2]{gal-19-wea}.
In the context of the heat equation ($F(u) = -\nabla u$), the numerical flux
\begin{equation}\label{eq:deffnumparab}
 F_{i+\frac 1 2}^n =\mathcal F_{i+\frac 1 2}^n(u_i^n,u_{i+1}^n) = -\frac {u_{i+1}^n - u_i^n} {h_{i+\frac 1 2}}\hbox{ with }h_{i+\frac 1 2} = x_{i+1}-x_{i},
\end{equation}
is consistent in the above sense, since, for a regular function $u$, 
\[
 \mathbb F(x_{i+\frac 1 2},t_n) = -\partial_x u(x_{i+\frac 1 2},t_n)
\]
and
\[
 \widetilde F_{i+\frac 1 2}^n(u(x_{i},t_n),u(x_{i+1},t_n)) = -\frac {u(x_{i+1},t_n) - u(x_{i},t_n)} {h_{i+\frac 1 2}} =  -\partial_x u(x_{i+\frac 1 2},t_n) - h_{i+\frac 1 2}\partial^2_{xx} u(c,t_n),\ c\in]x_i, x_{i+1}[.
\]
We notice that in this case the flux depends on the choice of the points $x_i$.

%

\section{Stability, compactness, convergence}

\subsection{Stability.} 
There are different notions of  stability of a numerical scheme.
The notion which is of interest in the context of the convergence of a numerical scheme for a general, possibly non linear, PDE  is an estimate on the approximate solutions, independently of the mesh. 
For instance, the $L^\infty$ stability of a linear FD scheme for a linear elliptic equation may be obtained by writing the scheme in matrix form and by obtaining a bound of the infinity norm of the inverse of this matrix.
Even though the FV approximate solutions are piecewise constant, the estimates are obtained in a norm which is in close relation with that which one uses for the estimates on the solutions of the continuous problem, and which depends of course on the considered problem;
the notion of stability in the FVM is therefore linked to the stability of the continuous problem.
Let us take two examples:

\begin{itemize}
 \item  \emph{the heat equation} on $[0,1]$ with  homogeneous Dirichlet conditions: 
the natural norm for the continuous problem is $L^2(H^1_0)$, and the associated discrete norm  corresponding to the choice \eqref{eq:deffnumparab} is $L^2(H^1_{0,d})$ with
\[
\Vert u \Vert_{H^1_{0,d}} = \Bigl(\sum_{i=0}^M h_{i+1/2} \Bigl( \dfrac {u_{i+1} - u_i} {h_{i+1/2}}\Bigr)^2 \Bigr)^{\frac 1 2},
\]
 where $M$ is the number of control volumes and $u_{0}=u_{M+1}=0$ and with $h_{i+1/2}= x_{i+1} - x_i$.
We denote by $u$ the weak solution of the heat equation and $u_{\mathrm{app}}$ the solution of the time-implicit scheme (implicit schemes are a natural choice for parabolic equations to avoid a condition of type $ \delta t \le C h^2 $).
The function $u_{\mathrm{app}}$ is a piecewise constant function which is equal to $u_i^n$ on the cuboid $]x_ {i-1/2}, x_{i + 1/2} [\times ] t_n, t_{n + 1} [$.
In the continuous PDE setting, an $ L^2 (H^1_0)$ estimate on $u$ is obtained by taking $u$ as a test function in the weak formulation of the heat equation, and integrating by parts.
Similarly, the $L^2(H^1_{0, d})$ estimate on $u_{\mathrm{app}}$ is obtained by multiplying the $i$-th discrete equation by $\delta t\, u_i^n$,  summing over $ i $ and  $ n $ and  performing discrete integrations by  parts (obtained by changes of indices in the sums).
An $L^2(L^2)$ estimate on the approximate solutions then follows with a discrete Poincaré inequality \cite[Lemma 9.1]{book}.

 \item \emph{the transport equation} on $\xR^n$, with initial condition $u_{\rm ini} \in L^\infty$:
 the natural norm for the continuous problem is $L^\infty$. 
 It is classical and easy to show that the scheme \eqref{vf}, with the upwind flux $F_{i+1/2} = a\, u_i$, is stable, see for example \cite[Lemma 20.1]{book}, under a CFL (for Courant Friedrichs Lewy, \cite{CFL-all}) condition.
 
 \end{itemize}

\subsection{The linear case: stability + conservativity +  flux consistency   $\implies$ convergence} 
Convergence for linear operators are often obtained through  error estimates (the same technique as in the proof of stability is applied to the error between the approximate solution and the exact one);
the compactness analysis of sequences of approximate solutions is another means, which also gives the existence of a solution (see \cite{CFL-all} for a seminal paper on this type of proof). 
Let us detail this second means, which extends to non linear operators.
Consider a sequence of approximate solutions, on meshes whose space and time steps tend to 0. 
The general principle of proof of convergence is then the following. 
Thanks to the stability, a uniform (with respect to the discretization parameters) estimate on the approximate solutions holds in a Lebesgue space, and so there exists a sub-sequence of this sequence which converges weakly (or $\star$-weakly) in that same Lebesgue space.
Each equation of the scheme is then multiplied by the interpolation of a regular test function and by the time step; the summation of all these equations over the time index and space indices, and an integration by parts (using the conservativity of the fluxes) are performed so as to shift the discrete derivatives from the discrete unknowns to the regular test function (the flux consistency which is used on the test functions is thus that given by the dual discrete  operator).
Since the problem is linear, the weak convergence of the approximate solutions suffices to pass to the limit in all terms; however the parabolic and hyperbolic cases exhibit some different difficulties.

\smallskip

\paragraph{- \emph{Elliptic or parabolic case}.} 
The passage to the limit can be performed thanks to the consistency of the flux, see e.g. \cite[Theorem 8.1]{book} for the 1D case. 
In the multidimensional case, the consistency of the flux for the heat equation (Laplace operator) is obtained for meshes which respect an orthogonality condition \cite[Definition 9.1]{book}. 
In this case, the resulting matrix is symmetric and the flux consistency is identical to that given by the dual discrete operator. 
Anisotropic operators and general meshes have been the object of several different works in the last decades, we refer to \cite{droniou2013review} for a review of these methods, several of them leading to non symmetric matrices. 
In this latter case, one can either prove convergence by an error estimate, using the flux consistency of the primal discrete operator, or a compactness method, using the flux consistency of the dual discrete operator, \ie on the test functions. 
In both cases the main difficulty is to establish the stability of the scheme.

\smallskip
 
\paragraph{- \emph{Hyperbolic case}. } 
For stability reasons, the numerical flux is upwinded and introduces an error term whose convergence to 0 must be shown. 
This term is a sum of  products of differences of the values of the solutions in neighbouring cells by discrete derivatives of the test function. 
On 1D or multiD Cartesian uniform or non uniform meshes, if $u_{\rm ini} \in L^\infty \cap BV$ (an additional argument is used to handle the case $u_{\rm ini} \in L^\infty$) this term is shown to tend to 0 thanks to a uniform $BV$ estimate on the approximate solutions (and on the continuous solution), see \cite{kuz-76-mon,CM} who consider the non linear case; this proof uses the ``TVD'' character (total variation diminishing) of the monotone schemes  \cite{harten}. 
Unfortunately, in the case of an unstructured mesh in multiD, even for a linear equation, even if $u_{\rm ini} \in L^\infty \cap BV$, it can be shown that the upwind scheme is not TVD; a counter-example is given in \cite{champier1992convergence}.
A ``weak BV'' estimate based on the numerical diffusion of the scheme is established therein in order to prove convergence.

\smallskip
 
Let us notice that in both the parabolic and hyperbolic cases, the Lax--Wendroff theorem is not used directly:
\begin{itemize}
 \item in the parabolic case, because the continuous flux function cannot be applied to the approximate solutions,
 \item in the hyperbolic case, because the approximate solutions only converge weakly, whereas the Lax--Wendroff theorem supposes a strong convergence. 
\end{itemize}
Nevertheless in both cases, the fundamental notions introduced in \cite{lax-60-sys} are used, namely: 
 
\begin{itemize}
  \item \emph{flux conservativity}: it is the property that leads to a weak form of the FV scheme by shifting the discrete derivatives of the discrete unknown to the discrete derivatives of the interpolant of the test function,
 \item \emph{flux consistency}: it is the property that is used to prove the fact that a limit of approximate solutions is a weak solution thanks to the convergence of the discrete derivatives of the interpolant of the test function to the exact derivatives of the same test function.
\end{itemize}
  
\subsection{The nonlinear case:  more compactness needed}
In the nonlinear case, weak convergence is not sufficient; indeed, if a sequence $(u_n)_\nnn$ converges only weakly in a Lebesgue space to a limit $u$, there is no reason why the sequence $(f(u_n))_\nnn$ should converge to $f(u)$, even weakly.
 
\smallskip

\paragraph{- \emph{Elliptic and parabolic equations}.}

In the elliptic or parabolic setting, one of the fundamental tools to obtain more compactness is Kolmogorov's compactness theorem, a consequence of which is that any bounded sequence of $L^p$, $1\le p <+\infty$ which is ``equicontinuous on average'' admits a convergent sub-sequence (see for example \cite[Theorem 8.16]{totoro}). 
Equicontinuity in mean amounts to showing that the difference between the function and its translates in time and space converges to $0$ in the $L^p$ norm, uniformly with respect to the time and space step, see for example \cite[Lemma 18.3]{book} in the case of a nonlinear parabolic equation of Stefan type. 
Once the compactness of the sequence of approximate solutions in $L^2(L^2)$ is proven, we can then exhibit a sub-sequence tending to $\bar u $ in $L^2(L^2)$.

By passing to the limit in the ``very weak'' form of the scheme (that is with the discrete divergence of the discrete normal derivatives of the test functions) we can then show, as in the linear case, that each term (in time and space) converges to the corresponding term in the ``very weak'' formulation of the continuous problem.  
In the case of the Stefan problem, namely $\partial_t u -\Delta \varphi(u)=0$, this gives the convergence of the approximate solutions to the exact solution if $\varphi$ is an increasing Lipschitz continuous fucntion. 
There is however an additional difficulty if $\varphi$ is only nondecreasing. 
In this case, it is possible to prove compactness in $L^2(L^2)$ of $\varphi(u_{\mathrm {app}})$ (where $u_{\mathrm{app}}$ is the approximate solution) but not of $u_{\mathrm{app}}$ for  which only an $L^2(L^2)$ bound holds.
It is possible however to conclude using the Minty trick (see \cite[Chapter 4]{book}).

\smallskip

 \paragraph{- \emph{Hyperbolic equations}}
 Now consider a nonlinear hyperbolic conservation law of the form
 \begin{align}
&\partial_t u (x,t) + \dive (f(u(x,t))) =0, \; x \in \xR^n, \\
&u(x,0) = u_{\rm ini}(x).
 \end{align}
If $f\in C^1(\xR, \xR^n)$ and $u_{\rm ini} \in L^\infty$, there exists a unique entropy solution of this problem \cite{krushkov}. 
In order to show that a scheme approximates this entropy solution, it is first shown that it satisfies a discrete entropy equation.
 
 The case of Cartesian meshes has been studied independently by Kuznetsov \cite{kuz-76-mon} and Crandall and Majda \cite{CM}.
 As in the linear case, if $u_{\rm ini}  \in BV$, a $BV$ estimate on the approximate solutions holds, uniformly with respect to the space and time step; Helly's lemma, which is itself a direct consequence of Kolmogorov's compactness theorem may then be invoked to obtain the convergence of a subsequence of the approximate solutions in $L^1(L^1)$, and one can then use the Lax--Wendroff theorem (which generalizes easily to the entropy formulation) \cite{CM} \cite[section 21.5]{book}).
It is also possible to handle the case $u_{\rm ini} \in L^\infty$, using a contraction principle in $L^1$ for the exact solution and for the approximate solution, see \cite{CM} for instance.
 
 In the general case of a non-Cartesian mesh, a suitable $BV$ estimate seems out of reach, and the proof of convergence is performed with the following steps.

 \begin{itemize}
 \item Consider a sequence of space discretizations, indexed by $m$, and of time steps $\delta t^{(m)}$~;  for any $m$, the mesh size is defined as the maximum diameter of the cells and denoted by $h^{(m)}$.
 Assume that $h^{(m)} \to 0$ and $\delta t^{(m)} \to 0$ as $m \to +\infty$. 
  \item Owing to the $L^\infty$ estimate on the sequence of approximate solutions $(u\exm)_\mnn$, there exists a sub-sequence which converges in a ``non linear weak sense'', \ie\ there exists $\mu \in L^\infty(\xR \times \xR_+ \times]0,1[)$ such that for any function $\psi \in C(\xR, \xR)$, 
  \[
    \int_\xR \int_{\xR_+} \psi(u\exm) \varphi \to \int_\xR \int_{\xR_+} \mu_\psi \varphi, \forall \varphi \in L^1(L^1) \mbox{ with } \mu_\psi = \int_0^1 \psi\bigl(\mu(x,t,\alpha)\bigr)\,\mathrm{d}\alpha.
  \]
 This notion of convergence is equivalent to the convergence to a Young measure \cite{DP}; it may seem a little simpler to handle in the sense that it involves a function, $\mu$, rather than a measure~; however this function depends on an additional parameter $\alpha \in ]0,1[$, which we will have to get rid of in order to reach  convergence to an entropy weak solution.
  \item Using the numerical diffusion of the scheme, we get a uniform weak $BV$ estimate on the sequence of approximate solutions; this estimate is called ``weak'' for two reasons: 
  on the one hand it involves the differences $|(f(u_K) - f(u_L) )\cdot \bfn_{KL}|$ and not the differences $|u_K - u_L|$, where $(K,L)$ denotes a pair of control volumes sharing a common interface~; on the other hand, it only requires that the sum of these differences does not blow up too fast: the difference between the discrete gradient of the interpolated test function and the gradient of the test function itself behaves like the mesh size, and to pass to the limit, we only need that the sum of the differences involving the discrete unknowns be bounded by a term in $C/ {h^{1-\varepsilon}}$ with $\varepsilon >0$ and $C>0$ independent of $h$.   
  
  \item Using the nonlinear weak convergence and the weak $BV$ estimate, we pass to the limit on the weak form of the discrete entropy and obtain a so-called ``process solution'' which is an entropy solution up to an integral with respect to the additional parameter $\alpha$.

  \item Starting from the discrete entropy inequalities that are verified by the scheme, the process solutions are shown to satisfy an entropy inequality.
A uniqueness result on the process solutions can then be obtained thanks to a variable doubling technique ``à la Kruskov''\cite{chinese}; this result differs from Di Perna's \cite{DP} in that it takes into account the initial condition in the weak entropy formulation, which allows to avoid the more restrictive conditions on the mesh \cite{coc-95-con}. 
  The uniqueness of the process solution entails that a process solution is the unique entropy weak solution.
It also yields the (strong) convergence, in $L^p$-spaces, of the approximate solution to the exact solution. 
\end{itemize}

 The proof of convergence of the FVM has been obtained for several other problems than the one considered here. 
 However, for systems of PDEs, it is often difficult to obtain compactness results, and LW--consistency then seems an interesting way to make sure that an eventual limit of the scheme is indeed a weak solution of the system.
 In the next section, we show how this is feasible even on staggered grids, which are often used in the numerical simulation of fluid flows.

\section{LW-consistency and staggered grids} \label{sec:MAC}

Rectangular staggered grids have been used since the sixties in fluid mechanics (in the wellknown MAC method \cite{har-65-num}) including environmental flows \cite{ara-81-pot}, see also \cite{pat-80-num}.
The mathematical analysis of the MAC scheme has been the object of several recent works, see e.g. \cite{gal-18-nsinc,gal-19-err}.
Systems of partial differential equations for which no existence or uniqueness is known have also been discretized on staggered grids. 
Let us mention in particular the compressible Euler equations:  
one of the advantages of using a staggered grid is to produce a scheme that can be mathematically proven to be asymptotically stable to the incompressible limit, as shown for the isentropic case in \cite{saleh}.
For such systems, a Lax--Wendroff type theorem is very useful; indeed, although true convergence cannot be proven for lack of compactness properties, such a theorem  allows to state that if the scheme converges, and provided some bounds on the approximate solutions are satisfied (these bounds are generally not attainable mathematically but can be verified numerically), the limit of the scheme is an entropy weak solution of the system, see \cite{her-15-swe}, \cite{her-21-cons}. 
The proof of this result may be obtained thanks to the generalization of the Lax--Wendroff theorem to general grids, which include staggered grids such as the MAC grid \cite{gal-21-wea}. 
One of the major additional difficulties for staggered grids is that the discrete unknowns are piecewise constant on different grids.

Consider a rectangular domain $\Omega \subset \xR^2$ (the 3D case can be tackled in the same way), and a possibly non-uniform rectangular grid. 
We denote by $\edgespart$ the set of edges of the mesh, and the internal edge separating the cells $K$ and $L$ is denoted by $\edge=K|L$ (see Figure \ref{fig:mac-space_disc}).
This mesh will be referred to in the following as the primal mesh, and denoted by $\mathcal P$.
Two dual meshes (three in 3D) are considered, each  consisting in a partition of $\Omega$ indexed by the vertical and horizontal elements of $\edgespart$, \ie\ $\Omega=\cup_{\edge \in \edgespart^{(i)}} D_\edge$, $i=1,2$, where $\edgespart^{(1)}$ (resp. $\edgespart^{(2)}$) denotes the set of vertical (resp. horizontal) edges.
The cells $(D_\edge)_{\edge \in \edgespart}$ are  referred to as the dual cells.
A half dual cell $D_{K,\edge}$ is  half of the rectangle $K$ with side $\edge$ (see Figure \ref{fig:mac-space_disc}). 
For an internal edge $\edge=K|L$, the dual cell $D_\edge$ is the subset of $K \cup L$ defined as $D_\edge=D_{K,\edge} \cup D_{L,\edge}$; for an external edge $\edge$ of a cell $K$, $D_\edge$ is the subset $D_{K,\edge}$ of $K$.

\begin{figure}[htb]
\begin{center}
\scalebox{0.9}{
%
\begin{tikzpicture} 
\draw[thick, color=black] (0.7, 0.5) -- (6.,0.5) ;
\draw[thick, color=black] (0.7, 2.) -- (6.,2.) ;
\draw[thick, color=black] (0.7, 4.5) -- (6.,4.5) ;
\draw[thick, color=black] (1.,0.2) -- (1.,4.8) ;
\draw[thick, color=black] (3.,0.2) -- (3., 4.8) ;
\draw[thick, color=black] (5.7,0.2) -- (5.7,4.8) ;
\node at (1.4,0.8){$\mathbf M$}; \node at (1.4,4.2){$\mathbf N$}; \node at (5.3,4.2){$\mathbf K$}; \node at (5.3,0.8){$\mathbf L$};
%
\draw[very thick, color=green!60!black] (3.,0.5) -- (3.,2.) node[midway, right]{$ \mathbf \edge = \mathbf M | \mathbf L$} ;
\draw[very thick, color=rougef] (3.,2.) -- (3.,4.5) node[midway, right]{$ \mathbf \edge' = \mathbf N | \mathbf K$} ;
\draw[very thick, color=bfonce] (1.,2.) -- (3.,2.) node[midway, above]{$ \mathbf \edgeb = \mathbf M | \mathbf N$} ;
\draw[very thick, color=orangec!80!black] (3.,2.) -- (5.7,2.) node[midway, above]{$ \mathbf \edgeb' = \mathbf K| \mathbf L$} ;
\end{tikzpicture}
%
\hspace{0.6cm} 
\begin{tikzpicture} 
\fill[color=rougef!20!white,opacity=0.3] (2.,2.) -- (2.,4.5) -- (4.35,4.5) -- (4.35,2.) -- (2.,2.);
\draw[very thick, color=rougef] (2.,2.) -- (2.,4.5) -- (4.35,4.5) -- (4.35,2.) -- (2.,2.);
\fill[color=green!30!black!20!white,opacity=0.3] (2.,0.5) -- (2.,2.) -- (4.35,2.) -- (4.35,0.5) -- (2.,0.5);
\draw[very thick,   color=green!20!black] (2.,0.5) -- (2.,2.) -- (4.35,2.) -- (4.35,0.5) -- (2.,0.5);
\node[color=green!30!black] at (3.9, 1.){$ {\mathbf D_\edge }$};
\node[color=orangec!80!black] at (3.9, 2.7){$ {\mathbf D_\edge' }$};

\draw[thick, color=black] (0.7, 0.5) -- (6.,0.5) ;
\draw[thick, color=black] (0.7, 2.) -- (6.,2.) ;
\draw[thick, color=black] (0.7, 4.5) -- (6.,4.5) ;
\draw[thick, color=black] (1.,0.2) -- (1.,4.8) ;
\draw[thick, color=black] (3.,0.2) -- (3., 4.8) ;
\draw[thick, color=black] (5.7,0.2) -- (5.7,4.8) ;
\node at (1.4,0.8){$\mathbf M$}; \node at (1.4,4.2){$\mathbf N$}; \node at (5.3,4.2){$\mathbf K$}; \node at (5.3,0.8){$\mathbf L$};
%
\draw[thick, color=rougef] (3.,2.) -- (3.,4.5) node[midway,right]{$\edge'$} ;
\draw[thick, color=green!60!black] (3.,0.5) -- (3.,2.) node[midway,right]{$\edge$} ;
\end{tikzpicture}
%
\hspace{0.6cm} 
\begin{tikzpicture} 
\fill[color=bclair,opacity=0.3] (1.,1.25) -- (1.,3.) -- (3.,3.) -- (3.,1.25) -- (1.,1.25);
\fill[color=orangec!30!white,opacity=0.3] (3.,1.25) -- (3.,3.) -- (5.7,3.) -- (5.7,1.25) -- (3.,1.25);
\draw[very thick,  color=bfonce]  (1.,1.25) -- (1.,3.) -- (3.,3.) -- (3.,1.25) -- (1.,1.25);
\draw[very thick,  color=orangec!80!black] (3.,1.25) -- (3.,3.) -- (5.7,3.) -- (5.7,1.25) -- (3.,1.25);
\node[color=bfonce] at (1.5, 2.6){$ {\mathbf D_\edgeb }$};
\node[color=rougef] at (5.2, 2.6){$ {\mathbf D_{\edgeb' }}$};

\draw[thick, color=black] (0.7, 0.5) -- (6.,0.5) ;
\draw[thick, color=black] (0.7, 2.) -- (6.,2.) ;
\draw[thick, color=black] (0.7, 4.5) -- (6.,4.5) ;
\draw[thick, color=black] (1.,0.2) -- (1.,4.8) ;
\draw[thick, color=black] (3.,0.2) -- (3., 4.8) ;
\draw[thick, color=black] (5.7,0.2) -- (5.7,4.8) ;
\node at (1.4,0.8){$\mathbf M$}; \node at (1.4,4.2){$\mathbf N$}; \node at (5.3,4.2){$\mathbf K$}; \node at (5.3,0.8){$\mathbf L$};
%
\draw[thick, color=bfonce] (1.,2.) -- (3.,2.) node[midway,above]{$\edgeb $} ;
\draw[thick, color=orangec!80!black] (3.,2.) -- (5.7,2.) node[midway,above]{$\edgeb'$} ;
\end{tikzpicture}
}
\end{center}
\caption{Primal and dual meshes and associated notations for the MAC case. \\- Left: the primal cells; the edges $\edge$ and $\edge'$ belong to $\edgespart^{(1)}$ and the edges $\edgeb$ and $\edgeb'$  to $\edgespart^{(2)}$. \\- Center: the dual cells associated to  $\edgespart^{(1)}$. \\- Right: the dual cells associated to  $\edgespart^{(2)}$.
}
\label{fig:mac-space_disc}
\end{figure} 
To illustrate the use of the generalized Lax--Wendroff theorem proven in \cite{gal-21-wea}, let us consider as a simple example the mass equation of, say, the compressible Euler equation:
\begin{equation} \label{def:mass}
   \partial_t\rho(\bfx,t) + \dive\bigl(\rho \bfu \bigr)(\bfx,t)=0 \quad(\bfx,t) \in \Omega\times ]0,T[,	 
\end{equation}
 where  $\partial_t \rho$ denotes the time derivative of the density $\rho$, and $\dive$ the space divergence.
The  scalar unknown $\rho$ is associated to the primal cells:
\begin{equation*}  
\rho(\bfx,t) = \rho_K^n \quad \mbox{for } \bfx \in K,\ K \in \mathcal P, \ t \in [t_n,t_{n+1}[,\ n \in \llbracket 0, N-1 \rrbracket.
\end{equation*}
The unknowns associated to the $i$-th component of $\bfu$ are located at the center of the edges of the $i$-th dual mesh.
The associated approximate vector function thus reads: $\bfu(\bfx,t) = (u_1(\bfx,t),\ u_2(\bfx,t))^t$  where, for $i=1,\ 2$,
\[  
u_i(\bfx,t) = u_\edge^n, \mbox{ for } \bfx \in D_\edge,\ \edge \in \edgespart^{(i)} \mbox{ and } t \in [t_n,t_{n+1}[,\ n \in \llbracket 0, N-1 \rrbracket.
\]
Let $\bfe^{(i)}$ denote the $i$-th unit vector; the discretization of \eqref{def:mass} reads:
\begin{multline*}
 \mathcal C(\rho,\bfu)_K^n = (\eth_t \rho)_K^n + \frac 1 {|K|} \sum_{\edge \in \edgespart(K)} |\edge|\ \bfF_\edge^n \cdot \bfn_{K,\edge}= 0,\mbox{ with } 
(\eth_t \rho)_K^n = \frac{\rho_K^{n+1}-\rho_K^n}{t_{n+1}-t_n}
\\
\mbox{ and } \bfF_\edge^n =  \rho_\edge^n \bfu_\edge^n,
\mbox{ where }\bfu_\edge^n  \mbox { is defined as } u_\edge^n  \ \bfe^{(i)} \mbox{ for } \edge \in \edgespart^{(i)},\ i = 1 \mbox{ or } 2,
\end{multline*}
and, for $\edge=K|L$, $\rho_\edge^n$ stands for a convex combination of $\rho_K^n$ and $\rho_L^n$ (for instance the upwind or a MUSCL choice with respect to $u_\edge^n$).
The initial value for the scalar unknown $\rho$ is defined by
\begin{equation}\label{eq:cond_ini}
\rho_K^0 = \frac 1 {|K|}\ \int_K \rho_0(\bfx) \dx. 
\end{equation}

\medskip
For $i = 1$, $2$, let $\bar h^{(i)}=\max \{ |\edge|,\ \edge \in \edgespart^{(i)} \}$ and $\underline h^{(i)}=\min \{ |\edge|,\ \edge \in \edgespart^{(i)} \}$.
We define the space step by $h(\mathcal P)=\max(\bar h^{(1)},\bar h^{(2)})$, and the time step by $\delta t = \max_{n \in \llbracket 0, N-1 \rrbracket} (t_{n+1}-t_n)$.
A sequence of grids is said quasi-uniform if the quotients $\bar h^{(1)}/\underline h^{(2)}$ and $\bar h^{(2)}/\underline h^{(1)}$ are bounded by a constant independent of the grid.

\begin{lemma}[Lax--Wendroff consistency for the mass equation, MAC grid] \label{lem:cons_stag}
Let a sequence of quasi-uniform MAC grids $(\mathcal P\exm)_\mnn$  and of time discretizations be given, with $h(\mathcal P\exm)$ and $\delta t\exm$ tending to zero;
let $(\rho\exm,\bfu\exm)_\mnn$ be the associated sequence of discrete functions.

We suppose that the sequences $(\rho\exm)_\mnn$ and $(\bfu\exm)_\mnn$ are bounded in $L^\infty(\Omega \times ]0,T[)$ and $L^\infty(\Omega \times ]0,T[)^2$ respectively, and that, when $m$ tends to $+\infty$,  they converge  in $L^p(\Omega \times ]0,T[)$ and $L^p(\Omega \times ]0,T[)^2$, $1 \leq p < +\infty$, to $\rho \in L^\infty(\Omega \times ]0,T[)$ and $\bfu \in L^\infty(\Omega \times ]0,T[)^2$ respectively.
Then $(\rho, \bfu)$ is a weak solution of \eqref{def:mass}, in the sense that, for any function $\varphi \in C^\infty_c(\Omega \times [0,T[)$,
\begin{multline*} 
\label{lw_ex}
\sum_{0\leq n \leq N-1} \delta t\exm \sum_{K \in \mathcal P\exm} |K| \, \mathcal C(\rho,\bfu)_K^n \, \varphi_K^n  \to -\int_\Omega \rho_0(\bfx)\ \varphi (\bfx,0) \dx
\\
- \int_0^T \int_\Omega \Bigl(\rho(\bfx,t)\ \partial_t \varphi(\bfx,t) + \ \bigl(\rho \bfu\bigr)(\bfx,t) \cdot \gradi \varphi(\bfx,t) \Bigr) \dx \dt
\quad \mbox{as } m \to + \infty,
\end{multline*}
where $\varphi_K^n$ stands for the mean value of $\varphi(\bfx,t_n)$ over the cell $K$, and so the right hand-side of this assertion vanishes.
\end{lemma}
The proof of Lemma \ref{lem:cons_stag} is given in \cite{gal-19-wea} and quite simple, especially using the {\em ad hoc} tools developed in \cite{gal-19-wea, gal-21-wea}.
However, the developed arguments may be extended to more complex operators, to deal for instance with the momentum and energy balance equations of the compressible Euler equations, even if the proofs are more tricky.

\section*{Conclusion}
 
In this paper, we have presented some concepts for the mathematical analysis of FV schemes, paying a special attention to the flux consistency issue and its most direct consequence, \ie\ the LW--consistency of the scheme.
We emphasize that we gave here a very partial picture of the mathematical world of finite volumes.
In many problems (some of them evoked here), the convergence of the discrete solution to a limit follows by compactness arguments in norms strong enough to "feed" the consistency study; {\em in fine}, this yields a stronger result than just the scheme consistency, namely the convergence (up to a subsequence if the uniqueness of the solution of the continuous problem is not known) of the numerical solutions to the (a) continuous one (see {\em e.g.} \cite{gal-18-nsinc}).
Many parabolic equations enter this framework, including, focusing on fluid flow simulations, incompressible, possibly variable density, or steady barotropic Navier-Stokes equations.
In some problems, the continuous solution may be reasonably supposed (or even proven) to be regular, and an error analysis is possible.

\medskip
The LW--consistency issue is especially important for practical applications in fluid flow simulations.
indeed, in many cases of interest, stronger results are out of reach, and this property is the only one left to mathematically support the design of schemes.
This is for instance the case for multi-dimensional flows governed by hyperbolic systems, as shallow-water equations, Euler equations or models for multi-phase flows.
For instance, the study evoked in Section \ref{sec:MAC} is motivated by such a situation: in the last ten years, a class of staggered schemes has been designed for hyperbolic flow problems \cite{her-15-swe,herbin2019decoupled,her-21-cons}, and implemented in the open-source software CALIF$^3$S developed at IRSN \cite{califs}; they are now routinely used for industrial safety applications as hydrogen explosion problems, supposing inviscid or at least vanishing viscosity flows.
The accuracy of the numerical schemes involved here is essentially supported by LW--consistency studies \cite{her-21-cons}.

\bibliographystyle{abbrv}
\bibliography{vf-cras}
\end{document}